\documentclass[12pt]{amsart}
\usepackage{amssymb}
\usepackage{amscd}

\usepackage[active]{srcltx}
\usepackage{epsfig}
\usepackage{epstopdf}
\usepackage{graphicx,color}

\nonstopmode

\DeclareMathOperator{\pnt}{\raise 0.5mm \hbox{\large\bf.}}

\newtheorem{theorem}{Theorem}[section]

\newtheorem{proposition}[theorem]{Proposition}
\newtheorem{corollary}[theorem]{Corollary}
\newtheorem{conjecture}[theorem]{Conjecture}

\theoremstyle{definition}
\newtheorem{remark}[theorem]{Remark}
\newtheorem{example}[theorem]{Example}

\begin{document}

\title[Pure $O$-sequences: known results, applications and open problems]{Pure $O$-sequences: known results,\\ applications and open problems}

\author[Juan Migliore]{Juan Migliore${}^*$}
\address{
Department of Mathematics \\
 University of Notre Dame \\
  Notre Dame,
IN 46556 \\
 USA}
\email{migliore.1@nd.edu}
\author[Uwe Nagel]{Uwe Nagel${}^{+}$}
\address{Department of Mathematics \\
University of Kentucky \\
715 Patterson Office Tower \\ Lexington, KY 40506-0027\\
 USA}
\email{uwe.nagel@uky.edu}
\author[Fabrizio Zanello]{Fabrizio Zanello}
\address{Department of Mathematical Sciences, Michigan Technological University, Houghton, MI 49931, USA}
\email{zanello@math.mit.edu}

\begin{abstract}
This note presents a discussion of the algebraic and combinatorial aspects of the theory of pure $O$-sequences. Various instances where pure $O$-sequences appear are described.  Several open problems  that deserve further investigation are also presented.
\end{abstract}

\thanks{${}^*$ The work for this paper was done while the first
author was sponsored by the National Security Agency under Grant
Number H98230-12-1-0204, and by the Simons Foundation under grant \#208579.\\
\indent ${}^+$ The work for this paper was done while the second author was
sponsored by the National Security Agency under Grant
Number H98230-12-1-0247, and by the Simons Foundation under grant \#208869.
}

\maketitle

\begin{center}
{\em Dedicated to David Eisenbud on the occasion of his 65th birthday}
\end{center}

\section{Introduction}
\label{sec-intro}

Pure $O$-sequences are fascinating objects that arise in several mathematical areas. They have been the subject of extensive research, yet our knowledge about pure $O$-sequences is limited. The goal of this note is to survey some of the known results and to motivate further investigations by pointing out connections to various interesting problems. Much of the material for this paper has been drawn from the recent monograph \cite{BMMNZ} by the authors with Mats Boij and Rosa Mir\'o-Roig, cited in this article  as BMMNZ. We often quote  results giving the name of the author(s) and the year the result was published, so as to also hint at the history of ideas in the development of the theory of pure $O$-sequences.

The multi-faceted interest in pure $O$-sequences is already indicated in their definition. On the one hand, a pure $O$-sequence can be defined as the vector whose entries record the number of monomials of a fixed degree in an order ideal generated by monomials of  the same degree. On the other, a pure $O$-sequence is the Hilbert function of a finite-dimensional graded algebra that is level, i.e., its socle is concentrated in one degree, and has monomial relations. There is an extensive literature on monomial ideals and an extensive literature on level algebras.  Pure $O$-sequences form a bridge between the two theories, and we will outline the work on pure $O$-sequences from this point of view.  But more than that, these sequences have a broad array of applications, and  occur in many settings, largely combinatorial ones.

In Section \ref{sec:level}, we review some of the basic results in the theory of pure $O$-sequences and focus on qualitative aspects of their shape. For instance, in some cases pure $O$-sequences are known to be unimodal; that is, they are first weakly increasing, and once the peak is reached they are weakly decreasing. However, pure $O$-sequences may fail to be unimodal,  even with arbitrarily many ``valleys." We include a discussion of recent results on conditions that force unimodality.

Connections to various combinatorial problems are the subject of Section \ref{sec:comb}. Face vectors of pure simplicial complexes are examples of pure $O$-sequences. In particular, the existence of certain block designs, such as Steiner systems,  is related to that of some pure $O$-sequences. As a special case, the existence of finite projective planes  is equivalent to the existence of particular pure $O$-sequences.

Another  challenging problem is Stanley's conjecture that the $h$-vector  of any matroid complex is a pure $O$-sequence. We discuss some recent progress. However, the conjecture remains open in general.

In Section \ref{sec:enum}, we describe results on the enumeration of pure $O$-sequences, with a focus on asymptotic properties. In particular, it follows that, when the number of variables is large, ``almost all'' pure $O$-sequences are  unimodal.

We conclude this note with a collection of open problems, most of which are mentioned in the earlier sections.


\section{Monomial  level algebras}
\label{sec:level}

A finite, nonempty set $X$ of (monic) monomials in the indeterminates $y_1,\dots,y_r$ is called a {\em monomial order ideal} if, whenever $M \in X$ and $N$ is a monomial dividing $M$, then $N \in X$.  The {\em h-vector} of $X$ is defined to be the vector $\underline{h} = (h_0 = 1, h_1,\dots,h_e)$ counting the number of monomials of $X$ in each degree.  A monomial order ideal, $X$, is called {\em pure} if all maximal monomials of $X$ (in the partial ordering given by divisibility) have the same degree.  A {\em pure $O$-sequence} is the $h$-vector of a pure monomial order ideal.  For reasons that will be clear shortly, we call $e$ the {\em socle degree} of $\underline{h}$.  The {\em type} of a pure $O$-sequence is the number of maximal monomials.

Notice that if we think of $y_1, \dots,y_r$ as the indeterminates of a polynomial ring $\mathcal R = K[y_1,\dots,y_r]$ over a field, the question of whether a given sequence is or is not a pure $O$-sequence does not depend on the choice of $K$. Thus, for many of our results it does not matter what we choose for $K$.  However, in some situations choosing a ``nice" field $K$ allows us to use special algebraic tools to say something about pure $O$-sequences, so in this case we make whatever additional assumptions we need for $K$.

We now let $R = K[x_1,\dots,x_r]$, where $K$ is an infinite field.  We will consider standard graded artinian $K$-algebras $A = R/I$, where $I$ will usually be a monomial ideal.  Without loss of generality we will assume that $I$ does not contain nonzero linear forms, so we will define $r$ to be the {\em codimension} of $A$.

Let $\mathcal R = K[y_1,\dots,y_r]$, and consider the action of $R$ on monomials of $\mathcal R$ by contraction.  By this we mean the action generated by
\[
x_i \circ y_1^{a_1} y_2^{a_2} \cdots y_r^{a_r} =
\left \{
\begin{array}{ll}
y_1^{a_1} y_2^{a_2} \cdots y_i^{a_1-1} \cdots y_r^{a_r}, & \hbox{if } a_i > 0, \\
0, & \hbox{if } a_i = 0.
\end{array}
\right.
\]

For a monomial ideal $I \subset R$, we define the {\em inverse system} to be the $R$-module  $I^\perp = \hbox{ann}_{\mathcal R}(I) \subset \mathcal R$.  One can check that $I^\perp$  consists of the monomials not in $I$ (identifying $x_i$ with $y_i$), and as such it can be viewed as a monomial order ideal.  Recalling that for a standard graded algebra $R/I$ the {\em Hilbert function} is defined to be $h_{R/I}(t) = \dim [R/I]_t$, we observe that the $h$-vector (as defined above) of the order ideal $I^\perp$ coincides with the Hilbert function of $R/I$.

Furthermore, $I^\perp$ is a pure monomial order ideal if and only if $R/I$ is a {\em level} algebra; that is, the socle of $R/I$ (i.e., the annihilator of the homogeneous maximal ideal of $R/I$) is concentrated in one degree, called the {\em socle degree} of $R/I$; it is necessarily the degree of the maximal monomials of $I^\perp$. The dimension of the socle as a $K$-vector space is equal to the type of the pure $O$-sequence.  See \cite{tony,IK} for more details on inverse systems.

Thus the study of pure $O$-sequences boils down to a study of the possible Hilbert functions of artinian monomial level algebras. In some cases we can rule out candidates for pure $O$-sequences by showing that there is not even a level algebra with that Hilbert function, but more often we need to use the structure of monomial algebras themselves.

The most basic tool is Macaulay's theorem to determine if the sequence is even an $O$-sequence; that is, to determine if it is the Hilbert function of {\em some} artinian algebra. We refer to  \cite{BH,mac} for details of Macaulay's theorem, but we recall the statement.  Let $n$ and $d$ be positive integers. There exist uniquely determined integers $k_d > k_{d-1} > \cdots > k_\delta \geq \delta \geq 1$ such that
\[
n = n_{(d)} = \binom{k_d}{d} + \binom{k_{d-1}}{d-1} + \cdots + \binom{k_\delta}{\delta}.
\]
This is called the {\em $d$-binomial expansion of $n$}.  We set
\[
({n_{(d)}})^1_1 = \binom{k_d+1}{d+1} + \binom{k_{d-1}+1}{d} + \cdots + \binom{k_\delta +1}{\delta +1}
\]
and $(0_{(d)})^1_1 = 0$, for each $d$. Then Macaulay's theorem is the following.

\begin{theorem}[Macaulay 1927] \label{macthm}
Let  $A$ be a standard graded algebra with Hilbert function $h_A(t) := h_t$.  Then for all $t\ge 1$, $h_{t+1} \leq ((h_t)_{(t)})^1_1$.
\end{theorem}

An {\em $O$-sequence} is a (possibly infinite) sequence of integers $(1,h_1,h_2,\dots )$  that satisfies the growth condition of Theorem \ref{macthm} for every value of $t$.  Thus the $O$-sequences are the sequences that occur as the Hilbert function of some standard graded algebra.

\begin{example}
The sequence $(1,3,6,8,8,10)$ is not a pure $O$-sequence because it is not even an $O$-sequence (the growth from degree 4 to degree 5 is too big).  Similarly, $(1,3,5,5,4,4)$ is an $O$-sequence but it is not a pure $O$-sequence because it is not the Hilbert function of a {\em level} algebra (see \cite{GHMS}).  Finally, $h=(1,3,6,10,15,21,28,27,27,28)$ is the Hilbert function of a level algebra, but it is not a pure $O$-sequence \cite{bernadette} because there is no {\em monomial} level algebra with this Hilbert function. In fact, $h$ has been the first  nonunimodal level Hilbert function discovered in codimension 3 (see the third author \cite{zanello1}), and Boyle \cite{bernadette} has shown that this is in fact the smallest possible such Hilbert function.
\end{example}

So the challenge is to determine what additional conditions on an $O$-sequence are imposed by requiring that it be the Hilbert function of an {\em artinian level monomial} algebra.  The first result, due originally to Stanley \cite{stanley1} with subsequent proofs given by J. Watanabe \cite{watanabe}, Ikeda \cite{ikeda},  Reid,  Roberts and  Roitman \cite{RRR}, Herzog and Popescu \cite{HP}, and  Lindsey \cite{lind}, concerns monomial complete intersections. (Note though that an equivalent property was proven earlier by de Bruijn, van Ebbenhorst Tengbergen and  Kruyswijk \cite{BTK}.)  This result requires us to introduce here the notion of the Weak and Strong Lefschetz Properties.  The consequence of this result for pure $O$-sequences of type 1 is perhaps not so critical, as alternative proofs could be given.  But its influence in the study of the SLP and the WLP in general, and on related topics in commutative algebra, can hardly be overstated. In the result below, and throughout the paper, we will call a sequence {\em unimodal} if it is nondecreasing up to some degree and then nonincreasing past that degree.  We will call it {\em strictly unimodal} if it is strictly increasing up to some degree, then possibly constant for some range, then strictly decreasing until it reaches zero, and then zero past that point. Unimodality is a central concept in combinatorics and combinatorial commutative algebra, but also in other branches of mathematics. See for instance the classical surveys of Stanley (\cite{St0}, 1989) and Brenti (\cite{Bre}, 1994).

\begin{theorem}[see above for sources]
\label{swrrr thm}

Let $R = K[x_1,\dots,x_r]$, where $K$ has  {\em characteristic zero}, and let $I$ be an artinian monomial complete intersection, i.e.,
\[
I = \langle x_1^{a_1},\dots,x_r^{a_r} \rangle.
\]
Let $L$ be a general linear form.  Then for any positive
integers $d$ and $i$, the homomorphism induced by multiplication by
$L^d$,
\[
\times L^d : [R/I]_i \rightarrow [R/I]_{i+d},
\]
has maximal rank.  (In particular, this is true when $d=1$.)  As a consequence, a pure $O$-sequence of type 1 is strictly unimodal.
\end{theorem}

When $I$ is an arbitrary artinian homogeneous ideal, the above maximal rank property for all $d$ and $i$ is called the {\em Strong Lefschetz Property (SLP)}, and the case $d=1$ is called the {\em Weak Lefschetz Property (WLP)}.  A consequence of the WLP is that in the range where $(\times L)$ is injective, we have a short exact sequence
\[
0 \rightarrow [R/I]_i \rightarrow [R/I]_{i+1} \rightarrow [R/(I,L)]_{i+1} \rightarrow 0.
\]
Thus, the {\em first difference} $\Delta h_{R/I} (t) = h_{R/I}(t) - h_{R/I}(t-1)$ is the Hilbert function of a standard graded algebra in this range, that is, it is again an $O$-sequence.  We then say that $h_{R/I}$ is a {\em differentiable $O$-sequence} in this range.  It also follows from this sequence that if the WLP holds, then  once the peak  of the Hilbert function is reached, the Hilbert function will be nonincreasing, and therefore the whole Hilbert function is unimodal.

\begin{remark}
  \label{rem:positve-char}
Simple examples show that Theorem \ref{swrrr thm} may fail in positive characteristic. It turns out that the question in which positive characteristics a monomial complete intersection has the weak or strong Lefschetz property leads to unexpected connections to the problem of determining the number of certain plane partitions,  lozenge tilings, or families of lattice paths (see \cite{CGJL, cook, CN, CN2, LZ}).
\end{remark}

One of the early important results on pure $O$-sequences is due to Hibi \cite{hibi}.

\begin{theorem}[Hibi 1989]
Let $\underline{h}$ be a pure $O$-sequence of socle degree $e$.  Then
\[
h_i \leq h_j
\]
whenever $0 \leq i \leq j \leq e-i$.  This has the following two important consequences:

\begin{itemize}
\item[(a)] $\underline{h}$ is {\em flawless}, i.e.,\ $h_i \leq h_{e-i}$ for all $0 \leq i \leq \lfloor \frac{e}{2} \rfloor$.

\item[(b)] The ``first half'' of $\underline{h}$ is nondecreasing:
\[
1 = h_0 \leq h_1 \leq h_2 \leq \cdots \leq h_{\lfloor \frac{e}{2} \rfloor}.
\]
\end{itemize}
\end{theorem}

This latter result was later improved by the following algebraic \emph{$g$-theorem} of Hausel \cite{hausel}:

\begin{theorem}[Hausel 2005] \label{hauselthm}
Let $A$ be a monomial Artinian level algebra of socle degree $e$. If the field $K$ has characteristic zero, then for a {\em general} linear form $L$, the induced multiplication
\[
\times L : A_j \rightarrow A_{j+1}
\]
is an injection, for all $j = 0,1,\dots,\lfloor \frac{e-1}{2} \rfloor$.  In particular, over any field, the sequence
\[
1, h_1 -1, h_2 - h_1 ,\dots, h_{\lfloor \frac{e-1}{2} \rfloor +1} - h_{\lfloor \frac{e-1}{2} \rfloor}
\]
is an $O$-sequence, i.e.,\ the ``first half'' of $\underline{h}$ is a {\em differentiable $O$-sequence}.
\end{theorem}

We have the following additional results on differentiability from  \cite{BMMNZ}:

\begin{theorem}[BMMNZ 2012]
\label{thm:diff-O}
\begin{itemize}

\item[(a)]
Every finite differentiable $O$-sequence $\underline{h}$  is the ``first half" of some pure $O$-sequence. (This is the converse of Hausel's theorem.)

\item[(b)]  In particular, any finite differentiable $O$-sequence is pure (by truncation).

\item[(c)] Any nondecreasing pure $O$-sequence of socle degree $\leq 3$ is differentiable.
\end{itemize}
\end{theorem}

It turns out that (c) is the best possible result in this direction:

\begin{proposition}[BMMNZ 2012]
There exist nondecreasing pure $O$-sequences of any socle degree $e \geq 4$ that are not differentiable.
\end{proposition}

\begin{example} \label{not diff}
We illustrate the preceding result with an example from \cite{BMMNZ}.  Observe first that the $h$-vector $\underline{h}' = (1,4,10,20,35)$ is a pure $O$-sequence since it is the $h$-vector of the truncation of a polynomial ring in four variables, $w,x,y,z$; the pure order ideal arises using all 35 monomials of degree 4 in $w,x,y,z$.  The $h$-vector $\underline{h}'' = (1,4,6,4,1)$ is also a pure $O$-sequence, since it is the  order ideal generated by a monomial $abcd$ in four new variables.  Now we work in a polynomial ring in the eight variables $w,x,y,z,a,b,c,d$, and we consider the pure order ideal generated by the above 36 monomials of degree 4.  The resulting $h$-vector $h$ is

\begin{center}
\begin{tabular}{cccccccccccc}
&1 & 4 & 10 & 20 & 35 \\
+&&4 & 6 & 4 & 1 \\ \hline
&1 & 8 & 16 & 24 & 36
\end{tabular}
\end{center}

\noindent Since the first difference of $h=(1,8,16,24,36)$ is $(1,7,8,8,12)$, which is not an $O$-sequence (because $12>(8_{(3)})^1_1=10$), we have constructed the desired example.
\end{example}

Putting aside the class of nondecreasing pure $O$-sequences, we now turn to the question of unimodality.  There are three factors that go into whether a nonunimodal example will exist: the codimension (i.e., the number of variables), the socle degree, and the type.  An easy application of Macaulay's theorem gives that any standard graded algebra of codimension two has unimodal Hilbert function, and in fact it is not hard to show that if $K[x,y]/I$ is level (monomial or not) then the Hilbert function is strictly unimodal (see e.g. Iarrobino \cite{Ia}, 1984). Hence the interesting questions arise for codimension $r \geq 3$.

We will begin with some results involving the socle degree (some of which will also bring in the codimension).   It follows from Hibi's result on flawlessness that any pure $O$-sequence of socle degree $\leq 3$ is unimodal.  The next case, socle degree 4, already is not necessarily unimodal, again thanks to an example from \cite{BMMNZ}:

\begin{example}
Observe that $\underline{h}' = (1,5,15,35,70)$ is a pure $O$-sequence, since it is the $h$-vector of the truncation of a polynomial ring in five variables, and as before $\underline{h}''=(1,4,6,4,1)$ is also pure, since it corresponds to the maximal monomial $abcd \in K[a,b,c,d]$. Hence, reasoning as above, we now consider one copy of $\underline{h}'$ and eleven copies of $\underline{h}''$ as $h$-vectors of pure $O$-sequences in twelve different rings, and we work in the tensor product of those rings. It follows that
\[
\underline{h}=(1,5,15,35,70)+ 11 \cdot (0, 4, 6,4,1)=(1,49,81,79,81)
\]
 is a nonunimodal pure $O$-sequence of socle degree 4.
 \end{example}

A natural question, then, is what is the smallest codimension for which nonunimodal pure $O$-sequences exist with socle degree 4. This remains open.  We do have the following results from \cite{bernadette}, however.

\begin{theorem}[Boyle 2012]
\begin{itemize}
\item[(a)] All pure $O$-sequences of socle degree $\leq 9$ in three variables are unimodal.

\item[(b)] All pure $O$-sequences of socle degree $\leq 4$ in four variables are unimodal.

\item[(c)] In four or more variables, there exist nonunimodal pure $O$-sequences in all socle degrees $\ge 7$.
\end{itemize}

\end{theorem}

\noindent In \cite{BZ}, Boij and the third author gave a nonunimodal pure $O$-sequence of codimension 3 and socle degree 12.  This is the smallest known example in codimension 3.  It follows from this and Boyle's result that in codimension 3, the only open cases are socle degrees 10 and 11.  Notice that, in codimension 4, the previous theorem leaves only open the socle degrees 5 and 6.

Now we turn to questions involving the type.  Of course it follows immediately from Theorem \ref{swrrr thm} that pure $O$-sequences of type 1 are unimodal.  What else can be deduced about pure $O$-sequences using the WLP?  A collection of results was obtained in \cite{BMMNZ} which showed, in some sense, the limits of the WLP in the study of pure $O$-sequences.

\begin{theorem}[BMMNZ 2012] \label{wlp situation}
Over a field of characteristic zero the following hold:
\begin{itemize}
\item[(a)] Any monomial artinian level algebra of type 2 in three variables has the WLP.  Thus a pure $O$-sequence of type 2 and codimension 3 is differentiable until it reaches its peak, is possibly constant, and then is nonincreasing until it reaches zero.

\item[(b)] Fix two positive integers $r$ and $d$.  Then all  monomial artinian level algebras of codimension $r$ and type $d$ possess the WLP if and only if at least one the following is true:

\begin{itemize}

\item[(i)] $r=1$ or 2;

\item[(ii)] $d=1$;

\item[(iii)] $r=3$ and $d=2$.

\end{itemize}

\end{itemize}
\end{theorem}

\noindent The proof of (a) was surprisingly long and intricate.  The main point of (b) is that in all other cases, we were able to show that artinian monomial level algebras exist that do {\em not} have the WLP.

Notice that Theorems \ref{swrrr thm}, \ref{hauselthm} and \ref{wlp situation}(a) require that $K$ have characteristic zero. The statements about injectivity and surjectivity require this property of the characteristic, and indeed a great deal of research has been carried out to see what happens when $K$ has positive characteristic; we refer to \cite{MN-exp} for an overview of these results.  The consequences on the shape of the pure $O$-sequences are indeed characteristic free, as has been noted above, since the Hilbert function of a monomial ideal does not depend on the characteristic.

If one is studying all artinian level monomial algebras of fixed type, Theorem \ref{wlp situation} is a serious limitation on the usefulness of the WLP.  However, in the study of pure $O$-sequences (e.g. to determine combinations of codimension, type and socle degree that force unimodality) it is still conceivable that the WLP will play a useful role.  As a trivial example, we know that monomial complete intersections in any codimension possess the WLP, thanks to Theorem \ref{swrrr thm}.  It is not known (except in codimension $\leq 3$, as noted below) whether {\em all} complete intersections have this property, but the knowledge in the special case of monomial ideals is enough to say that all complete intersection Hilbert functions are unimodal.  Perhaps a similar phenomenon will allow the WLP to continue to play a role in the study of pure $O$-sequences. A first approach using this philosophy was obtained by Cook and the second author \cite{CN-lifting}, where they {\em lifted} a monomial ideal to an ideal of a reduced set of points in one more variable, showed that the general artinian reduction has the WLP, and concluded that the Hilbert function of the original monomial algebra is unimodal, regardless of whether it has the WLP or not.

For this reason we mention a useful tool in studying the WLP, that was introduced by the first and second authors with T. Harima and J. Watanabe in \cite{HMNW} and whose study was continued by H. Brenner and A. Kaid in \cite{BK}.  This is the study of the {\em syzygy bundle}, and the use of the Grauert-M\"ulich theorem (see \cite{OSS}).  It was used in \cite{HMNW} to show that {\em any} complete intersection $I$ in $K[x,y,z]$ over a field of characteristic zero has the WLP.  The idea is to restrict to a general line, say one defined by a general linear form $L$.  The key is that the restricted ideal $(I,L)/(L)$ (which now has codimension 2, hence has a Hilbert-Burch matrix) should have minimal syzygies in two consecutive degrees, at most.  The idea of \cite{HMNW} was that for height 3 complete intersections,  this information can be obtained by considering the module of syzygies of $I$, sheafifying it to obtain the syzygy bundle, and applying the Grauert-M\"ulich theorem to the general line defined by $L$.  Of course this introduces questions about the semistability of the syzygy bundle, which we mostly omit here.  In a more general setting, the following result from \cite{BK} summarizes the idea nicely, at least for codimension three.

\begin{theorem}[Brenner-Kaid  2007]
Let $I = \langle f_1,\dots,f_k \rangle \subset K[x,y,z] =R$ be an artinian homogeneous ideal whose  syzygy bundle $\mathcal S$ is semistable on $\mathbb P^2$.  Then
\begin{itemize}
\item[(a)] If the restriction of $\mathcal S$ splits on a general line $L$ as
\[
\mathcal S_L \cong \bigoplus_{i=1}^s \mathcal O_L (a+1) \oplus \bigoplus_{i=s+1}^{k-1} \mathcal O_L(a),
\]
then $A = R/I$ has the WLP.

\item[(b)] If the restriction of $\mathcal S$ splits on a general line $L$ as
\[
\mathcal S_L \cong \mathcal O_L(a_1) \oplus \cdots \oplus \mathcal O_L(a_{n-1})
\]
with $a_1 \geq a_2 \geq \dots \geq a_{n-1}$ and $a_1 - a_{n-1} \geq 2$, then $A = R/I$ does not have the WLP.
\end{itemize}

\end{theorem}

Other applications of this approach, for higher codimension, can be found for instance in \cite{MMN2}.

Returning to unimodality questions, tools other than the WLP will also be needed.  We have the following theorems of Boyle \cite{bernadette}, which relied on  decomposition results and an analysis of complete intersection Hilbert functions but did not use the WLP.

\begin{theorem}[Boyle 2012]
\begin{itemize}
\item[(a)] In codimension 3, all pure $O$-sequences of type 3 are strictly unimodal.

\item[(b)] In codimension 4, all pure $O$-sequences of type 2 are strictly unimodal.

\end{itemize}
\end{theorem}

\noindent A natural question is whether all pure $O$-sequences of type 2 and arbitrary codimension are unimodal.  Also, for any fixed codimension, it is an interesting problem to determine which types force unimodality. The ``record" for the smallest known nonunimodal example in codimension 3 is type 14, given in \cite{BMMNZ}.

Finally, one can ask ``how nonunimodal''  a pure $O$-sequence can be.  The answer is ``as nonunimodal as you want.'' We have:

\begin{theorem}[BMMNZ 2012] \label{arb peaks}
For any integers $M\ge 2$ and $r \geq 3$, there exists a pure $O$-sequence in $r$ variables which is nonunimodal and has exactly $M$ maxima.
\end{theorem}

Of course the ``price" in Theorem \ref{arb peaks} is paid in having a large socle degree and type.

In fact, even \emph{Cohen-Macaulay $f$-vectors} (i.e., the face vectors of Cohen-Macaulay  simplicial complexes, which are a much smaller subset of pure $O$-sequences) can be nonunimodal with arbitrarily many peaks (see \cite{fz}).  This result considerably extends Theorem \ref{arb peaks}, even though, unlike for arbitrary pure $O$-sequences, here the number of variables becomes necessarily very large as the number of peaks increases.
\smallskip

A good topic to build a bridge between the algebraic and the combinatorial sides of the theory of pure $O$-sequences is the \emph{Interval Conjecture}.

The \emph{Interval Property} (IP) was  introduced in 2009 by the third author \cite{Za}, where he conjectured its existence for the set of Hilbert functions of level --- and, in a suitably symmetric way, Gorenstein --- algebras. Namely, the IP says that if two (not necessarily finite) sequences, $h$ and $h'$, of a  class $S$ of integer sequences coincide in all entries but one, say
$$h=(h_0,\dots ,h_{i-1},h_i,h_{i+1},\dots ) {\ }{\ }\text{and}{\ }{\ }h'=(h_0,\dots ,h_{i-1},h_i+\alpha ,h_{i+1},\dots ),$$
for some index $i$ and some positive integer $\alpha $, then the sequences
$$(h_0,\dots ,h_{i-1},h_i+\beta ,h_{i+1},\dots )$$
are also in $S$, for all  $\beta =1,2,\dots , \alpha -1.$

Given that level and Gorenstein Hilbert functions are nearly impossible to characterize, the IP appears to be both a very natural property and one of the strongest structural results  that we might hope to achieve for the set of such sequences.

For example, it is proved in  \cite{Za} that the IP holds for all Gorenstein Hilbert functions of socle degree 4. Since Gorenstein Hilbert functions are symmetric,  this means that, for any fixed $r$, $(1,r,a,r,1)$ is  Gorenstein if and only if $a$ ranges between some minimum possible value, say $f(r)$, and $\binom{r+1}{2}$. This latter is  the maximum allowed by a polynomial ring in $r$ variables, and is  achieved by the so-called \emph{compressed} Gorenstein   algebras (see, e.g., the 1984 papers of Fr\"oberg-Laksov and Iarrobino  \cite{FL,Ia} or the third author's works \cite{Za00,Za000}).

Notice that the existence of the IP for these Hilbert functions is especially helpful in view of the fact that, for most codimensions $r$, the  value of $f(r)$ is not known. In fact, such Gorenstein Hilbert functions are ``highly'' nonunimodal. Asymptotically, we have $$\lim_{r\rightarrow \infty} \frac{f(r)}{ r^{2/3}}= 6^{2/3},$$ as proved by these three authors in \cite{MNZ1} (2008), solving a longstanding conjecture of Stanley \cite{St3} (see also \cite{MNZ3},  2009, for some broad generalizations).

The Interval Property is still wide open today for both level and Gorenstein algebras.

In a more  combinatorial direction, in BMMNZ 2012 the IP has then also been conjectured for: 1) pure $O$-sequences (under the name ``ICP''); and 2) the $f$-vectors of pure simplicial complexes, a topic of discussion of the next section.

As for pure $O$-sequences,  the ICP has been proved in a number of special cases. Most importantly, it is known when the socle degree is at most 3, in any number of variables.

\begin{theorem}[BMMNZ 2012]
The ICP holds for the set of all pure $O$-sequences of socle degree $e\le 3$.
\end{theorem}

Thanks to  this  result,  H.T. H\`a, E. Stokes and the third author \cite{HSZ} have recently developed a new approach leading to a proof of Stanley's matroid $h$-vector conjecture in Krull-dimension 3, as we will see in the next section.

While the ICP remains open in most instances --- e.g., in three variables --- it must be pointed out that, just recently, it has been \emph{disproved} in the four variable case by Constantinescu and Varbaro (see \cite[Remark 1.10]{CV}), who found the following counterexample.

\begin{example}[Constantinescu-Varbaro 2012]  \label{CV ex}
Consider the pure order ideals generated by $\{x^3 y^2 z, x^3 y t^2, x^3 z^2 t\}$ and $\{x^4 y^2, x^3 y z t, x^2 z^2 t^2\}$. Their $h$-vectors are the pure $O$-sequences $(1, 4, 10, 13, 12, 9, 3)$ and $(1, 4, 10, 13, 14, 9, 3)$. However, an exhaustive computer search over all sets of three monomials of degree 6 in four variables reveals that the sequence $(1, 4, 10, 13, 13, 9, 3)$ is not pure, contrary to the  ICP.

 It is worth remarking that $h = (1,4,10,13,13,9,3)$ is, however, a level $h$-vector, and so this does not provide a counterexample to the IP for arbitrary level algebras.  Indeed, $h$ is   the $h$-vector of a level algebra in four variables whose inverse system is generated by two sums of sixth powers of six general linear forms each, and the sixth power of one general linear form.
\end{example}

As for pure $f$-vectors, the IP is still wide open, and little progress has been made so far.

In general, at this time  it is still unclear what the exact scope  of the Interval Property is, and if it can also be of use  in other areas of combinatorial algebra or even enumerative combinatorics. It is well known to hold, e.g., for the set of Hilbert functions of graded algebras of any Krull-dimension (see Macaulay's theorem),  the $f$-vectors of arbitrary simplicial complexes (the Kruskal-Katona theorem), and the $f$-vectors of Cohen-Macaulay complexes (BMMNZ 2012). Instead, the IP fails quite dramatically, for example, for matroid $h$-vectors, which are  conjecturally another subset of pure $O$-sequences, as we will see in the next  section (we refer to BMMNZ 2012 and \cite{St3} for details). Stanley and the third author \cite{SSZZ} recently looked at the IP in the context of \emph{$r$-differential posets}, a class of ranked posets generalizing the Young lattice of integer partitions and the Young-Fibonacci lattice. Here, even though the IP  fails in general, it might be a reasonable property to conjecture, for instance, for $r=1$, which is the most natural class of differential posets.


\section{Pure $O$-sequences and combinatorics}
\label{sec:comb}

Much of the motivation for the study of pure $O$-sequences comes from combinatorics, and in this section we give an overview of this side of the theory. In order to put in context the definition of a pure $O$-sequence given in the introduction, we quickly recall the notion of posets and order ideals. For an introduction to this theory, we refer to Chapter 3 of Stanley's new edition of ``EC1'' (\cite{EC1}, 2012). A \emph{poset} (short for \emph{partially ordered set}) is a set $S$ equipped with a binary relation, ``$\le ,$'' that is: 1) reflexive (i.e., $a\le a$ for all $a\in S$); 2) antisymmetric ($a\le b\le a$ implies $a=b$); and 3) transitive ($a\le b\le c$ implies $a\le c$).

An \emph{order ideal} in a poset $S$ is a subset  $I$ of $S$ that is closed with respect to ``$\le .$'' That is, if $t\in I$ and $s\le t$, then $s\in I$. Thus, our monomial order ideals are  the (finite) order ideals  of the poset $P$ of all monomials in the polynomial ring $R = K[y_1,\dots,y_r]$, where the binary relation of $P$ is divisibility. Notice that $P$ is a \emph{ranked} poset, where the rank of a monomial is its degree in $R$. Therefore, Macaulay's $O$-sequences  are exactly the possible rank functions of the  order ideals of $P$, since every Hilbert function satisfying Macaulay's theorem can be achieved by a monomial algebra. Similarly, as we have seen, a {pure} $O$-sequence is the rank function of some monomial order ideal whose  \emph{generators} (i.e., the antichain of maximal monomials) are all of the same degree.

Another fundamental class of order ideals are those contained in the \emph{Boolean algebra} $B_r$, the poset of all subsets of $\{1,2,\dots,r\}$, ordered by inclusion. By identifying the integer $i$ with a vertex $v_i$, the order ideals of $B_r$ are usually called \emph{simplicial complexes} (on $r$ vertices).

The elements of a simplicial complex $\Delta$ are dubbed \emph{faces}, and the maximal faces are the \emph{facets} of $\Delta$. The \emph{dimension of a face} is its cardinality minus 1, and the \emph{dimension of $\Delta$} is the largest of the dimensions of its faces.

Notice that if we  identify $v_i$ with a variable $y_i\in R=K[y_1,\dots,y_r]$, then simplicial complexes also coincide with the order ideals of $P$ generated by \emph{squarefree} monomials. In particular, if we define as \emph{pure} those simplicial  complexes whose facets have all the same dimension, then clearly their rank vectors, called \emph{pure $f$-vectors}, are the special subset of pure $O$-sequences that can be generated by squarefree monomials.

\begin{example}
The simplicial complex
$$\Delta =\left\{ \{1,2,3\},\{2,3,4\},\{1,2\},\{1,3\},\{2,3\},\{2,4\},\{3,4\},\{1\},\{2\},\{3\}, \{4\},\emptyset  \right\}$$
is the order ideal of $B_4$ generated by $\{1,2,3\}$ and $\{2,3,4\}$.

Thus, $\Delta$ is a pure complex of dimension 2, whose pure $f$-vector is $f_{\Delta}=(1,4,5,2)$.

Equivalently, $f_{\Delta}$ is  the pure $O$-sequence generated by the two squarefree monomials $y_1y_2y_3$ and $y_2y_3y_4$.
\end{example}

Similarly to Macaulay's theorem for arbitrary $O$-sequences,  we know a characterization of the class of pure $f$-vectors thanks to the classical \emph{Kruskal-Katona theorem} (see e.g. \cite{St3}). However,  analogously, things become dramatically more complicated (hopeless, we should say) when it comes to attempting a characterization of \emph{pure} $f$-vectors.

In the last section of  BMMNZ 2012, we have  begun a study of pure $f$-vectors, but still very little  is known today beyond what is known for arbitrary pure $O$-sequences.

Besides their obvious intrinsic importance --- simplicial complexes are a central object in algebraic combinatorics, combinatorial algebra and topology, just to name a few subjects --- pure $f$-vectors also carry fascinating  applications. It is on their connections to finite geometries and design theory  that we want to focus in the next portion of this section.

It will follow from our discussion, as probably first observed by Bj\"orner (\cite{bjorner}, 1994), that a characterization of pure simplicial complexes and their $f$-vectors would imply, for instance, that of all Steiner systems, and as a further special case, a classification of all finite projective planes, one of the major open problems in geometry.
\smallskip

A \emph{Steiner system} $S(l,m,r)$ is an $r$-element set $V$, together with a collection of $m$-subsets of $V$, called \emph{blocks}, such that every $l$-subset of $V$ is contained in exactly one block. Steiner systems are a special family of the so-called \emph{block designs}. We refer our reader to the two texts \cite{codi,liro}, where she can find a truly vast amount of information on  combinatorial designs.
For instance, a Steiner \emph{triple} system (STS) is a Steiner system $S(2,3,r)$, while $S(3,4,r)$ is dubbed a Steiner \emph{quadruple} system, where $r$ is the \emph{order} of the system.

Since we are dealing with maximal sets of the same cardinality ($m$, in this case), it is clear that if we identify each element of $V$ with a variable $y_i$, then the existence of Steiner systems (and similarly for other block designs) will be equivalent to the existence of certain pure $f$-vectors.

\begin{example}
Let us consider  STS's of order 7, i.e., $S(2,3,7)$. Constructing such a design is tantamount to determining a family of squarefree degree 3 monomials  of $R=K[y_1,y_2,\dots,y_7]$, say $M_1,M_2,\dots,M_t$,  such that each squarefree degree 2 monomial of $R$ divides exactly one of the $M_i$.

Clearly, since there are $\binom{7}{2}=21$ squarefree degree 2 monomials in $R$, if the $M_i$ exist, then $t=21/\binom{3}{2}=7$. In other words, $S(2,3,7)$ exists if and only if
$$f=(1,7,21,7)$$
is  a pure $f$-vector.

Notice also that $f$ exists as a pure $f$-vector if and only if it exists as a pure $O$-sequence, since for seven degree 3 monomials to have a total of 21 degree 2 divisors, each needs to have exactly three linear divisors, i.e.,  it must be squarefree.

It is easy to see that an STS of order 7, and so the pure $f$-vector $f=(1,7,21,7)$, do indeed exist, using the monomial order ideal generated by:
$$y_1y_2y_3,{\ } y_3y_4y_5,{\ }y_3y_6y_7,{\ }y_1y_4y_7,{\ }y_2y_4y_6,{\ }y_2y_5y_7,{\ }y_1y_5y_6.$$
\end{example}

Some simple numerical observations show that a necessary condition for an STS of order $r$ to exist is that $r$ be congruent to 1 or  3 modulo 6, and it is a classical result of Kirkman (\cite{kirk}, 1847) that this is also sufficient. In other words,
$$f=\left(1,r,\binom{r}{2}, \binom{r}{2}/3 \right)$$
is a pure $f$-vector, if and only if it is a pure $O$-sequence, if and only if $r$ is congruent to 1 or  3 modulo 6.

A different, and more challenging, problem is the classification of all Steiner systems, up to isomorphism. Even the existence of  particular systems sometimes brings into the story a nontrivial amount of interesting algebra. As an illustration, we mention here the case of the Steiner systems $S(4,5,11)$, $S(5,6,12)$, $S(3,6,22)$, $S(4,7,23)$, and $S(5,8,24)$, which are intimately connected to the first  sporadic finite simple groups ever discovered, called the \emph{Mathieu groups} (see E.L. Mathieu \cite{Mat1,Mat2}, 1861 and 1873). These five groups --- denoted respectively by $M_{11}$,  $M_{12}$,  $M_{22}$,  $M_{23}$, and  $M_{24}$ ---  in fact arise as the \emph{automorphism groups} of the above Steiner systems (i.e., the transformations of the systems that preserve the blocks).

There exists only one STS of order 7, which is called the \emph{Fano plane} (see the figure below) for reasons that will be clear in a minute. In other words, the  seven monomials of the previous example are, up to isomorphism, the only possible set of generators for a pure order ideal in $K[y_1,y_2,\dots,y_7]$ with $(1,7,21,7)$ as its pure $O$-sequence.

\begin{center}

\includegraphics[width=2in]{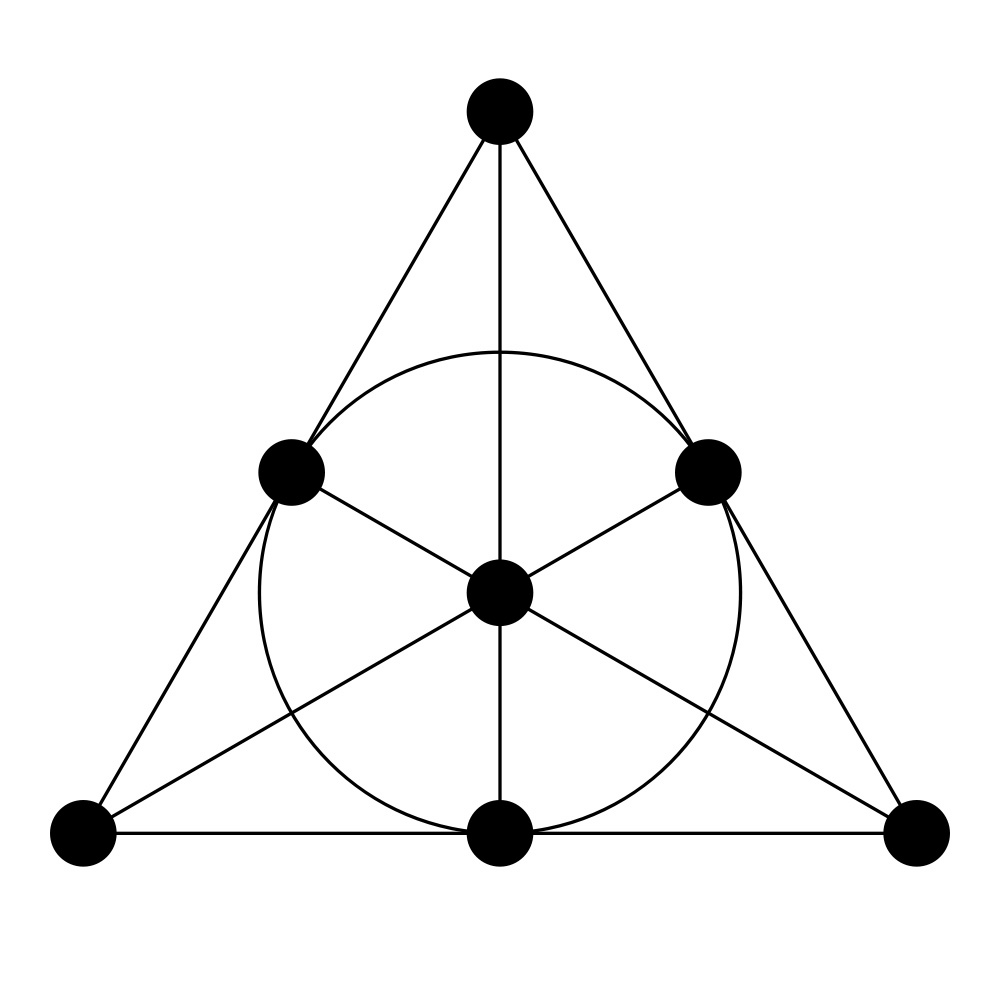}

\end{center}

\begin{picture}(150,0)
\put(223,148){$y_1$}
\put(164,95){$y_2$}
\put(256,95){$y_5$}
\put(226,72){$y_4$}
\put(134,37){$y_3$}
\put(209,25){$y_7$}
\put(286,37){$y_6$}

\end{picture}

Also for $r=9$, there exists a unique STS. However,  it is reasonable to believe that the  number of nonisomorphic STS increases extremely quickly for $r$ large. For instance, there are 80  nonisomorphic STS of order 15, and there are $11,\hspace{-.05cm}084,\hspace{-.05cm}874,\hspace{-.05cm}829$ of order 19.

Similarly, the possible orders of Steiner quadruple systems are known and nicely characterized, since the obvious necessary conditions  again turn out to be also sufficient: $S(3,4,r)$ exists if and only if its order $r$ is congruent to 2 or  4 modulo 6, as proved by Hanani (\cite{hana}, 1960). In other words, reasoning as above, since there are four possible 3-subsets of any given 4-set, we have that
$$f=\left(1,r,\binom{r}{2}, \binom{r}{3}, \binom{r}{3}/4\right)$$
is a pure $f$-vector if and only if it is a pure $O$-sequence, if and only if $r$ is congruent to 2 or  4 modulo 6.

However, as we increase the cardinality $m$ of the blocks, things become more and more obscure. This is due to the high complexity of computing combinatorial designs over a large vertex set, as well as to the lack of a general theory.

Already  for Steiner quintuple systems, the trivial necessary conditions ($r$ congruent to 3 or 5, but not to 4, modulo 6) are no longer sufficient. For example, no Steiner quintuple system $S(4,5,17)$ exists (see \cite{ospo}, 2008). The smallest value of $r$ for which the existence of $S(4,5,r)$ is currently open is $21$. In other words, it is unknown whether
$$\left(1,21,\binom{21}{2},\binom{21}{3},\binom{21}{4},\binom{21}{4}/5\right)=\left(1,21,210,1330,5985,1197\right)$$
is a pure $O$-sequence.

Perhaps the best-known family of examples of Steiner systems is that of {finite projective planes}, so they deserve a special mention here. Recall that a \emph{projective plane} is a collection of points and lines such that any two lines ``intersect at'' exactly one point, and any two points ``lie on'' exactly one line (one also assumes that there exist four points no three of which are collinear, in order to avoid uninteresting pathological situations).

If the projective plane is {finite}, it can easily be seen that the number of points is equal to the number of lines, and that this number must be of the form $q^2+q+1$, where the integer $q$ is the \emph{order} of the plane. Further, in a projective plane of order $q$, any line contains exactly $q+1$ points, and by duality,  any point is at the intersection of exactly $q+1$ lines.  In other words, finite projective planes  are the Steiner systems  $S(2,q+1,q^2+q+1)$. The reader may want to consult, e.g.,  \cite{De} for an introduction to this area.

Thus, the above example of a Steiner system $S(2,3,7)$, the Fano plane, is the unique smallest possible projective plane.

Similarly to how we argued earlier in terms of pure $f$-vectors, one  can show that a projective plane of order $q$ exists if and only if
$$h= \left(1,q^2+q+1,(q^2+q+1)\binom{q+1}{2},(q^2+q+1)\binom{q+1}{3},\dots ,(q^2+q+1)\binom{q+1}{q+1} \right) $$
is a pure $f$-vector, if and only if it is a pure $O$-sequence.

A major open problem in  geometry asks for a classification of all finite projective planes, or even just of the possible values that $q$ may assume. Conjecturally, $q$ is always the power of a prime, and it is a standard  algebraic  exercise, using finite field theory, to construct a projective plane of any order $q=p^n$.

The best general necessary condition known today on $q$ is still the following  theorem from \cite{BR,CR}:

\begin{theorem}[Bruck-Ryser-Chowla \cite{BR,CR}]
If $q$ is the order of a projective plane and $q$ is congruent to 1 or 2 modulo 4, then $q$ is the sum of two squares.
\end{theorem}

Thus, for instance, as a consequence of the Bruck-Ryser-Chowla Theorem, no projective plane of order 6 exists. In other words,
$$(1,43, 903, 1505, 1505, 903, 301, 43)$$
is not a pure $O$-sequence. However, already ruling out the existence of projective planes of order 10 has required a major computational effort (see Lam \cite{Lam}, 1991). The case $q=12$ is still open.

Notice that, at least for certain values of $q$, the number of nonisomorphic projective planes of order $q$ can be very large, and a general classification seems entirely out of reach. The smallest  $q$ for which there exists more than one nonisomorphic projective plane is  $9$, where the four possible cases  were already known to O. Veblen (\cite{VW}, 1907).
\smallskip

The second important application of pure $O$-sequences that we want to discuss brings our attention to a very special class of simplicial complexes, called \emph{matroid} complexes. Matroids are ubiquitous in mathematics, where they often show up in surprising ways (see \cite{NN,Ox,Wh,Wh2}).

The algebraic theory of matroids began in the same 1977 seminal paper of Stanley \cite{St1} that introduced pure $O$-sequences.  A finite matroid can be naturally identified with a pure simplicial complex over $V=\{1,2,\dots ,r\}$, such that its restriction to any subset of $V$ is also a pure complex.

One associates, to any given simplicial complex $\Delta$ over $V=\{1,2,\dots,r\}$, the following squarefree monomial ideal in $S =K[x_1, \dotsc, x_r]$, where $K$ is a field:
$$I_\Delta = \left \langle x_F = \prod_{i \in F} x_i \mid F \not\in \Delta \right \rangle.$$
\noindent
$I_\Delta$ is called the \emph{Stanley-Reisner ideal} of the complex $\Delta$, and the quotient algebra $S/I_\Delta$ is its \emph{Stanley-Reisner ring}.

It is a standard fact of combinatorial commutative algebra (see e.g. \cite{St3}) that the Stanley-Reisner ring $S/I_\Delta$  of a matroid complex is Cohen-Macaulay and  level, although of course of positive Krull-dimension (except in degenerate cases). Thus, the $h$-vector of $S/I_\Delta$ is  level, since it is the $h$-vector of an artinian reduction of $S/I_\Delta$. However, even though $S/I_\Delta$ is presented by monomials, notice that its artinian reductions will in general be far from  monomial, for they require  taking quotients by  ``general enough'' linear forms.

The following spectacularly simple conjecture of Stanley (\cite{St1}, 1977) predicts that, for any matroid complex $\Delta$, we can nonetheless find \emph{some} artinian {monomial} level algebra having the $h$-vector of $S/I_\Delta$ as its $h$-vector:

\begin{conjecture}\label{stst}
Any matroid $h$-vector is a pure $O$-sequence.
\end{conjecture}

The problem of characterizing matroid $h$-vectors  appears to be once again hopeless, and Conjecture \ref{stst} has motivated much of the algebraic work done on matroids over the past 35 years  (see, as a highly nonexhaustive list, \cite{Ch,Ch2,CV,DKK,HaS,Me,Ohh,Pr,Sch,Sto,Sw,Sw2}).

The main approach to Conjecture \ref{stst} has  been, given the $h$-vectors of a certain class of matroids,  to explicitly produce some pure monomial order ideals  having those matroid $h$-vectors as their pure $O$-sequences. Recently,  H.T. H\`a, E. Stokes and the third author \cite{HSZ} introduced a ``more abstract'' approach to Conjecture \ref{stst}. Their main idea, inspired by the latest progress on pure $O$-sequences made in BMMNZ 2012, and in particular the proof of the Interval Conjecture (ICP) in socle degree 3, has been to try to reduce Stanley's conjecture, as much as possible, to one on the properties of pure $O$-sequences, thus avoiding  explicit construction of a monomial ideal for each matroid $h$-vector.

The  approach of \cite{HSZ} has already led  to a proof of Conjecture \ref{stst} for all matroid complexes of Krull-dimension at most 2 (the dimension 1 case, that had been the focus of a large portion of the thesis \cite{Sto}, simply followed in  a few lines).

\begin{theorem}[H\`a-Stokes-Zanello \cite{HSZ}]
All matroid $h$-vectors $(1,h_1,h_2,h_3)$ are pure $O$-sequences.
\end{theorem}

More generally,  the following is a first concrete, if still  tentative, general approach to Conjecture \ref{stst} (see \cite{HSZ}).

Assuming  Conjecture \ref{stst} holds for all matroid complexes whose deletions with respect to any vertex are cones (which may not be too difficult to show with the techniques of paper \cite{HSZ}), Conjecture \ref{stst} is true in general under the following two  natural (but still too bold?)  assumptions:
\begin{enumerate}
\item[(A)] \emph{Any matroid $h$-vector is {differentiable} for as long as it is nondecreasing.} (In fact, incidentally, would a \emph{$g$-element}, that Hausel \cite{hausel} and Swartz \cite{Sw} proved to exist in the ``first half'' of a matroid, carry on all the way?)
\item[(B)] \emph{Suppose that the \emph{shifted sum}, $h''=(1,h_1+1,h_2+h_1',\dotsc,h_e+h_{e-1}')$, of two pure $O$-sequences $h$ and $h'$ is differentiable for as long as it is nondecreasing. Then $h''$ is also a pure $O$-sequence.}
\end{enumerate}


\section{Enumerations of pure $O$-sequences}
\label{sec:enum}

As we have seen above, pure $O$-sequences arise in several areas, yet their properties are not well understood, and there are other important questions that should be addressed even if a classification is not available. For example, one would like to estimate the number of pure $O$-sequences of given codimension and socle degree. What happens asymptotically? Moreover, we have seen that pure $O$-sequences can be as far from being unimodal as we want. Nevertheless, one may ask: What are the odds for a pure $O$-sequence to be unimodal?

In order to discuss such questions, let us denote by $O(r, e)$, $P(r, e)$, and $D(r, e)$ the sets of $O$-sequences, pure $O$-sequences, and differentiable $O$-sequences, respectively,  that have codimension $r$ and socle degree $e$. Recall that given two  functions $f, g: {\mathbb R} \to {\mathbb R}$,  one says that $f$ is {\em asymptotic} to $g$, and writes $f(r)\sim_r g(r)$, if $\lim_{r}f(r)/g(r)=1$. All limits are taken for $r$ approaching infinity.

Consider now an $O$-sequence $(1, r-1, h_2,\ldots,h_e)$ in $O(r-1, e)$. Integrating it, that is, passing to $(1, r, r-1 + h_2,\ldots,h_{e-1} + h_e)$,   provides  a differentiable $O$-sequence in $D(r, e)$. Thus, since finite differentiable $O$-sequences are pure by Theorem \ref{thm:diff-O}(b), we have the following inclusions:
\[
O(r-1, e) \hookrightarrow  D(r, e) \subset P(r, e) \subset O(r, e).
\]
Results by Linusson (see \cite{Li})  imply that, for  $r$ large, the cardinalities of $O(r-1, e)$ and $O(r, e)$ are asymptotically equal. It follows that in large codimensions almost all $O$-sequences are pure. More precisely, one has:

\begin{theorem}[BMMNZ 2012]
\label{almost} Fix a positive integer $e$.  Then, for  $r$ large, almost all $O$-sequences of socle degree $e$ are differentiable. Namely,
\[
\# O(r, e) \sim_r \# P(r, e) \sim_r \# D(r, e)  \sim_r c_e \cdot r^{\binom{e+1}{2}-1},
\]
\label{c_e}
where
\[
c_e= \frac{{\displaystyle \prod_{i=0}^{e-2}\binom{\binom{e+1}{2}-\binom{i+1}{2}-1}{i}}}{{\displaystyle \left ( \binom{e+1}{2}-1 \right ) !}}.
\]
\end{theorem}

Since pure $O$-sequences are Hilbert functions of level algebras, we immediately get the following consequence.

\begin{corollary}[BMMNZ 2012]
\label{level}
Fix a positive integer $e$. Let $L(r,e)$ be the set of  {\em level} Hilbert functions of codimension $r$ and socle degree $e$. Then, for $r$ large, almost all level sequences are pure and unimodal, and
\[
\#L(r,e) \sim_r c_e r^{\binom{e+1}{2}-1}.
\]
\end{corollary}

The outcome changes drastically if we fix as additional parameter the socle type $t$.
Since each monomial of degree $t$ is divisible by at most $t$ distinct variables, we observe:

\begin{proposition}[BMMNZ 2012]
\label{000}
Let $P(r,e,t)$ be the set of pure $O$-sequences of codimension $r$, socle degree $e$, and type $t$. Then $\# P(r,e,t)=0$ for $r>te$,  this bound being sharp.
\end{proposition}

Note that, in contrast to this result,  there is no analogous restriction on level Hilbert functions. For example, Gorenstein algebras have type 1 and admit any positive socle degree and codimension. In fact, asymptotically, the number of their Hilbert functions is known.

Recall that an \emph{SI-sequence} of socle degree~$e$ is an $O$-sequence  that is symmetric about $\frac{e}{2}$ and differentiable up to degree $\lfloor \frac{e}{2} \rfloor$. Initially, Stanley and Iarrobino (see \cite{St3}) had hoped that all Hilbert functions of  Gorenstein algebras were SI-sequences. Although this is not true (see \cite{stanley1} for the first counterexample), it is {\em almost true}! In fact, any differentiable $O$-sequence in $D(r, \lfloor \frac{e}{2} \rfloor)$ can be extended to a symmetric sequence, so that the result is an SI-sequence. Moreover, every SI-sequence is the Hilbert function of some Gorenstein algebra (see \cite{CI, harima, MN3}). Obviously, the first half of the Hilbert function of a Gorenstein algebra is an $O$-sequence. Taken together, it follows that the number of Gorenstein Hilbert functions that are not SI-sequences is negligible:

\begin{theorem}[BMMNZ 2012]
\label{gorsi} Fix a positive integer $e$. Let $G(r,e)$ be the set of Gorenstein Hilbert functions of codimension $r$ and socle degree $e$, and let $SI(r,e)$ be the set of SI-sequences of codimension $r$ and socle degree $e$. Then, for $r$ large, almost all  Gorenstein Hilbert functions are SI-sequences. More precisely,
$$\# G(r,e)\sim_r \# SI(r,e)\sim_r c_{\lfloor e/2 \rfloor}r^{\binom{{\lfloor e/2 \rfloor}+1}{2}-1}.$$
\end{theorem}

Returning to pure  $O$-sequences, it would be very interesting to determine, or at least to find a good estimate of, the number of pure $O$-sequences of codimension $r$, socle degree $e$, and type $t$. This seems a  difficult problem. However, in the simplest case, where $t = 1$, there is an easy combinatorial answer. In fact,  any  pure $O$-sequence of type 1 is the Hilbert function of a (complete intersection) algebra whose inverse system is a monomial of the form $y_1^{a_1} \cdots y_r^{a_r}$, where $a_1+ \dots + a_r=e$ and $a_i\geq 1$ for all $i$. We may assume that $a_1 \ge a_2 \ge \cdots \ge a_r$, so that $(a_1,\ldots,a_r)$ is a \emph{partition} of $e$. Since it is easy to see that distinct partitions lead to different Hilbert functions, we arrive at the following result.

\begin{proposition}[BMMNZ 2012]
\label{p_r} $\# P(r,e,1)=p_r(e)$, the number of partitions of the integer $e$ having exactly $r$ parts.
\end{proposition}

We now consider a slightly different asymptotic enumeration question. Fix positive integers $e$ and $t$. Since the number of monomials dividing $t$ monomials of degree $e$ is finite, there are only finitely many pure $O$-sequences of socle degree $e$ and type $t$. Denote their set by $P(\ast ,e,t)$. Determining $\# P(\ast ,e,t)$ exactly seems out of reach. However, one may hope to  at least find its order for  $t$ large. The first interesting case, namely $e = 3$, has been settled.

\begin{theorem}[BMMNZ 2012]
\label{thmAsymptoticSocleDegree3}
Let $\# P(\ast ,3,t)$ denote the number of pure $O$-sequences of socle degree 3 and  type $t$.  Then
\[
\lim_{t \rightarrow \infty} \frac{\# P(\ast ,3,t)}{t^2} = \frac{9}{2}.
\]
\end{theorem}

Its proof gives some further information. Consider a pure $O$-sequence $(1,r,a,t)$.  Then
\[
r \leq a \leq 3t,
\]
by Hibi's theorem and the fact that $t$ monomials of degree 3 are divisible by at most $3t$ distinct quadratic monomials.
It follows that $\# P(\ast ,3,t) \le \frac{9}{2} t^2$. To see that the two functions are in fact asymptotically equal, notice that, for fixed $t$, the possible values of $r$ and $a$ fall into one of the following three regions, illustrated in the figure below:

{Region I}: \quad $t \leq r \leq a \leq 3t$;

{Region II}: \quad  $0 < r < t \leq a \leq 3t$;

{Region III}: \quad  $0 < r \leq a < t$.

\

\bigskip

\begin{center}
\begin{picture}(170,160)(10,10)
\thicklines
\put (10,10){\vector(0,1){170}}
\put (10,10){\vector(1,0){170}}
\put (10,10){\line(1,1){160}}
\put (180,0){$r$}
\put (0,180){$a$}
\thinlines
\put (60,10){\line(0,1){150}}
\put (160,10){\line(0,1){150}}
\put (10,60){\line(1,0){150}}
\thicklines
\put (10,160){\line(1,0){150}}
\thinlines
\put (60,-2){$t$}
\put (156,-2){$3t$}
\put (-4,156){$3t$}
\put (0,60){$t$}
\put (87,120){{\Large I}}
\put (27,120){{\Large II}}
\put (15,40){{\Large III}}
\end{picture}
\end{center}
\bigskip  \bigskip

Using superscripts to denote the sets of pure $O$-sequences in each region, it is shown in BMMNZ  that

\begin{equation*}\label{222}
\lim_{t \rightarrow \infty} \frac{\# P^{(I)}(\ast ,3,t)}{t^2} = 2, \quad \lim_{t \rightarrow \infty} \frac{\# P^{(II)}(\ast ,3,t)}{t^2} = 2, \quad \text{and } \lim_{t \rightarrow \infty} \frac{\# P^{(III)}(\ast ,3,t)}{t^2} = \frac{1}{2}.
\end{equation*}

The arguments use in a crucial way the interval property for pure $O$-sequences of socle degree 3. Unfortunately, this property fails in general. Nevertheless it would be very interesting to extend the above results to socle degree $e \ge 4$.


\section{Open problems}

In this section we collect a few interesting problems that remain open in the area of pure $O$-sequences.  Most of them have been discussed in the previous sections, but some are related problems that were not addressed above.

\begin{enumerate}
\item What is the largest type $t$ for which all pure $O$-sequences  are unimodal (independently of the codimension or socle degree)?  Even proving that $t \ge 2$, i.e., that pure $O$-sequences of type 2 are unimodal in any codimension, would be very interesting.

\item For a fixed codimension $r$, what is the largest type $t$ for which all pure $O$-sequences are unimodal?

\item For a fixed codimension $r$, what is the largest socle degree $e$ for which all pure $O$-sequences are unimodal?

\item What is the smallest codimension $r$ for which there exists a nonunimodal pure $O$-sequence of socle degree 4?

\item Determine asymptotically the number $\# P(\ast ,e,t)$ of pure $O$-sequences of socle degree $e$ and type $t$, for  $t$ large. What is the order of magnitude of $\# P(\ast ,e,t)$?

\item The first example of a nonunimodal pure $O$-sequence (due to Stanley \cite{St1}, 1977) was $(1,n=505,2065, 3395, 3325, 3493)$, which is in fact the $f$-vector of a Cohen-Macaulay simplicial complex, hence in particular a pure $f$-vector. What is the smallest number of variables $n$ (i.e., the number of vertices of the  complex) allowing the existence of a nonunimodal pure $f$-vector?

A. Tahat \cite{tahat} has recently determined  the sharp lower bound $n=328$ for  nonunimodal \emph{Cohen-Macaulay} $f$-vectors of socle degree 5 (the socle degree of  Stanley's original example), and has produced examples in larger socle degree with $n$ as low as 39.

\item Stanley's Twenty-Fifth Problem for the  year 2000 IMU Volume ``Mathematics: frontiers and perspectives''  \cite{st-imu} asks, among a few other things: are all matroid $f$-vectors unimodal, or even {log-concave}? What about matroid $h$-vectors?

Notice  that  matroid $f$-vectors  are a much smaller subset of Cohen-Macaulay $f$-vectors. For a major recent breakthrough on this problem, see  Huh  \cite{Hu}, and Huh and Katz \cite{HK}; as a consequence of their work, log-concavity (hence  unimodality) is now known for all $f$-vectors (see Lenz \cite{lenz}, 2012) and $h$-vectors (Huh \cite{huh3}) of  \emph{representable} matroids. Notice also that, in general, the log-concavity of matroid $h$-vectors would imply the log-concavity of matroid $f$-vectors, as proved by Dawson \cite{jer} (1984) and Lenz \cite{lenz}.

\end{enumerate}

\section*{Acknowledgements} We wish to thank David Cook II and Richard Stanley for helpful comments. The third author also thanks J\"urgen Bierbrauer for an interesting discussion on the connections between group theory and Steiner systems.  We thank the referee for a careful reading of the paper.


\end{document}